\documentclass[12pt,bezier]{article}
\usepackage{times}
\usepackage{booktabs}
\usepackage{pifont}
\usepackage{floatrow}
\usepackage{caption}
\usepackage{mathrsfs}
\usepackage[fleqn]{amsmath}
\usepackage{amsfonts,amsthm,amssymb,mathrsfs,bbding}
\usepackage{txfonts}
\usepackage{graphics,multicol}
\usepackage{graphicx}
\usepackage{color}
\usepackage{longtable}
\usepackage[colorlinks,
            linkcolor=black,
            anchorcolor=black,
            citecolor=black]{hyperref}
\usepackage{caption}
\usepackage{cite}
\usepackage{latexsym,bm}
\usepackage{indentfirst}
\usepackage{color}
\usepackage[T1]{fontenc}
\usepackage[utf8]{inputenc}
\usepackage{authblk}
\evensidemargin=\oddsidemargin\topmargin=-1.5cm

\newtheorem{lemma}{Lemma}[section]
\newtheorem{theorem}{Theorem}[section]
\newtheorem{cor}{Corollary}[section]

\newtheorem{conj}{Conjecture}

\theoremstyle{definition}

\addtocounter{section}{0}

\begin{document}

\title{\bf\Large On the sum of the first two largest signless Laplacian eigenvalues of a graph\footnote{ Corresponding author Email:\ cxhe@usst.edu.cn (C.-X. He).\\
		This research is supported by  National Natural Science
		Foundation of China (No.12271362).   }}
\author[1]{Zi-Ming Zhou}
\author[1]{Chang-Xiang He}
\author[2]{Hai-Ying Shan}
\affil[1]{College of Science, University of Shanghai for Science and Technology, Shanghai, P. R. China.}
\affil[2]{School of   Mathematical Science, Tongji University, Shanghai, P. R. China.}
 
\renewcommand\Authands{ and }

\date{}

\maketitle

\begin{abstract}
For a graph $G$, let $S_2(G)$ be the sum of the first two largest signless Laplacian eigenvalues of $G$,    
and $f(G)=e(G)+3-S_2(G)$.  Oliveira, Lima, Rama and Carvalho conjectured that $K^+_{1,n-1}$ (the star graph with an additional edge)  is the unique graph with minimum value of $f(G)$ on $n$ vertices. In this paper, we prove this conjecture,    which also confirm  a conjecture for the upper bound of $S_2(G)$ proposed by Ashraf et al.  
\medskip

\begin{flushleft}
\end{flushleft}

\end{abstract}

\section{Introduction}

Graphs in this paper are connected and simple. Let $G=(V(G),E(G))$ be a   graph, the number of its vertices and edges are denoted by $n(G)$ (or $n$ for short) and $e(G)$, respectively. Let $A(G)$ be the adjacency matrix of $G$ and $D(G)=diag(d(v_1),d(v_2),\ldots,d(v_n))$ be the diagonal matrix of vertex degrees of $G$.
The Laplacian matrix and the signless Laplacian matrix of $G$ are defined as $L(G)=D(G)-A(G)$ and $Q(G)=D(G)+A(G)$, respectively.
It is well known that both $L(G)$ and $Q(G)$ are positive semidefinite and hence their eigenvalues are nonnegative real numbers.
The eigenvalues of $L(G)$ and $Q(G)$ are called  Laplacian eigenvalues and  signless Laplacian eigenvalues of $G$, respectively,
and are denoted by $\mu_1(G)\geq\mu_2(G)\geq\cdots\geq\mu_n(G)$
(or $\mu_1\geq\mu_2\geq\cdots\geq\mu_n$ for short)
and $q_1(G)\geq q_2(G)\geq\cdots\geq q_n(G)$
(or $q_1\geq q_2\geq\cdots\geq q_n$ for short), respectively.

In \cite{Brouwer}, Brouwer proposed the following conjecture for Laplacian eigenvalues.

\begin{conj}\label{Brouwer Conjecture}( \cite{Brouwer})
\emph{For any graph $G$ with $n$ vertices and for any $k\in\{1,2,\ldots,n\}$,
\begin{equation*}
\sum_{i=1}^k\mu_i(G)\leq e(G)+\binom{k+1}{2}.
\end{equation*}
}
\end{conj}

It has been proved that Conjecture \ref{Brouwer Conjecture} is true for several classes of graphs   such as trees \cite{Haemers}, threshold graphs \cite{Haemers}, unicyclic graphs \cite{DuZhou}, \cite{WangHuang}, bicyclic graphs \cite{DuZhou}, regular graphs \cite{Mayank} and split graphs \cite{Mayank}. Recently, Cooper \cite{Cooper} showed that Conjecture \ref{Brouwer Conjecture} is true for planar graphs when $k\geq 11$ and for bipartite graphs when $k \geq \sqrt{32}n$.
For $k = 2$, Haemers et al.\cite{Haemers} showed that Conjecture \ref{Brouwer Conjecture} is true.
In \cite{LG}, Li and Guo proposed the full Brouwer's Laplacian spectrum conjecture and proved the conjecture holds for $k=2$ with the unique extremal graph.

Let $S_k(G)$
be the sum of the first $k$ largest signless Laplacian eigenvalues of $G$.  Ashraf et al.\cite{Ashraf} proposed the following conjecture.

\begin{conj}\label{Ashraf Conjecture}(\cite{Ashraf})
\emph{For any graph $G$ with $n$ vertices and for any $k\in\{1,2,\ldots,n\}$,
\begin{equation*}
S_k(G)\leq e(G)+\binom{k+1}{2}.
\end{equation*}
}
\end{conj}

For bipartite graphs, 
  Conjecture \ref{Ashraf Conjecture}   is equivalent to  Conjecture \ref{Brouwer Conjecture}  since $L(G)$ and $Q(G)$ are similar \cite{DCve}.
Ashraf et al.\cite{Ashraf} proved that Conjecture \ref{Ashraf Conjecture} is true for regular graphs, they also proved that Conjecture \ref{Ashraf Conjecture} is true for $k=2$, 
but the key lemma (\cite{Ashraf}, Lemma 13) they used is incorrect which has a counterexample in\cite{Zheng}.
Zheng\cite{Zheng} proved that Conjecture \ref{Ashraf Conjecture} is true for all connected triangle-free graphs when $k = 2$.
Therefore, for general graphs Conjecture \ref{Ashraf Conjecture} is still open when $k=2$.

For a graph $G$, define  
\begin{equation*}
f(G)=e(G)+3-S_2(G).
\end{equation*} 
Oliveira et al.\cite{Oliveira} conjectured that $K^+_{1,n-1}$ (the star graph with an additional edge) minimizes $f(G)$ among all graphs $G$ on $n$ vertices. 

\begin{conj}(\cite{Oliveira})\label{Oliveira}
For any graph $G$ on $n\geq 9$ vertices. Then \begin{equation*}
f(G)\geq f(K^+_{1,n-1}).
\end{equation*} Equality if and only if $G\cong K^+_{1,n-1}$.
\end{conj}

In this paper, we give a proof  of Conjecture \ref{Oliveira}.

\begin{theorem}\label{them1}
For any graph $G$ with $n$ vertices, $f(G)\geq f(K^+_{1,n-1})$ with equality if and only if $G\cong K^+_{1,n-1}$.
\end{theorem}
 By computing the value of $f(K^+_{1,n-1})$, one can obtain a tight upper bound of  $S_2(G)$ which confirm Conjecture \ref{Ashraf Conjecture} in the case  $k=2$.



The organization of the remaining of the paper is as follows.
In Section 2, we give some lemmas which are useful in the proof of Theorem \ref{them1}. In Section 3, we give the proof of Theorem \ref{them1}.

\section{Preliminaries}

We begin this section with some definitions and notation. For vertex $v\in V(G)$,  the set of its  neighbors  in  $X\subset V(G)$,  denoted by  $N_X(v)$, for convenience, abbreviate $N_{V(G)}(v)$ as $N(v)$. For two graphs $G_1$ and $G_2$, the union of $G_1$ and $G_2$, denoted by $G_1\cup G_2$, is the graph with vertex
set $V(G_1)\cup V(G_2)$ and edge set $E(G_1)\cup E(G_2)$.
If $H$ is a subgraph of $G$, let $G-E(H)$ denote the graph obtained  from $G$ by removing the edges in $H$.

 For any square matrix $M$, we use $P(M,x)$ to denote the characteristic polynomial of $M$. For graph $G$, we write the characteristic polynomial of its signless Laplacian matrix $P(Q(G),x)$ as $P_Q(G,x)$.

{\it $H$-join operation} of graphs: Let $H$ be a graph with $V(H)=\{1,2,...,k\}$, 
and $G_i$ be disjoint $r_i$-regular graphs of order $n_i\ (i = 1,2,...,k)$. 
Then the graph $\bigvee_H\{G_1,G_2,...,G_k\}$ is formed by taking the graphs $G_1,G_2,...,G_k$ and joining every vertex of $G_i$ to every vertex of $G_j$ whenever $i$ is adjacent to $j$ in $H$.    
Define $N_i=\sum_{j\in N_H(i)}n_j$ for each $i\in V(H)$ and
\begin{align*}
(Q^\pi)_{ij}=
\begin{cases}
2r_i+N_i, &i=j,\\
n_j, &ij\in E(H), \\
0, &otherwise.
\end{cases}
\end{align*}

\begin{lemma}[\cite{WH}]\label{join}
Let $H$ be a graph with $V(H)=\{1,2,...,k\}$, and $G_i$ be a $r_i$-regular graph of order $n_i\ (i = 1,2,...,k)$.
If $G=\bigvee_H\{G_1,G_2,...,G_k\}$, then
\begin{align*} 
P_Q(G,x)=P(Q^\pi,x)\prod_{i=1}^k\frac{P_Q(G_i,x-N_i)}{x-2r_i-N_i}.
\end{align*}
\end{lemma}

{\it Lexicographically ordered} $k$-sets: Fix a number $k$ between 1 and $n$, inclusive. Consider all possible sets of $k$ distinct numbers between 1 and $n$. Order the elements of each set by  using a {\it lexicographic} or {\it dictionary } order. Denote
these sets by $S_1, S_2, \ldots, S_{\tbinom{n}{k}}$ so that  $S_1< S_2< \ldots< S_{\tbinom{n}{k}}$.


{\it Compound matrix}: Let $A_{i,j}$ be the $k \times k$ matrix formed by the intersection, in $A$, of the rows whose numbers are in the set  $S_i$ and the columns whose numbers are in the set $S_j$. Let $A$ be a complex matrix  of order $n$, define $C_k(A)$ as  its $k$-th  {\it compound matrix}, which is the matrix of size $ \tbinom{n}{k} \times \tbinom{n}{k}$ whose $(i,j)$-th entry is the determinant of the matrix $A_{i,j}$.

{\it Additive compound matrix}: The matrix $C_k(I+tA)$ is an $ \tbinom{n}{k} \times \tbinom{n}{k}$ matrix,
each entry of which is a polynomial in $t$ of degree at most $k$.
The $k$-th  {\it additive compound matrix} is defined as
\begin{align*} 
\Delta_k(A)= \frac{d}{dt}C_k(I_n+tA)\big\vert_{t=0}.
\end{align*}

For additive compound matrix, the following lemma is a well-known result.

\begin{lemma}[\cite{MF}]\label{MF}
If $\lambda_1,\ \lambda_2,\ldots, \lambda_n$ are all eigenvalues of $A$,
then the eigenvalues of $\Delta_k(A)$ are the $\tbinom{n}{k}$  distinct sums of the $\lambda_i$ taken $k$ at a time.
\end{lemma}

\begin{cor}\label{c2}
If $q_1,\ q_2,\ldots, q_n$ are signless eigenvalues of $G$,
and $P(Q^\pi,x)$ contains the largest two roots of $P_Q(G,x)$,
then the largest eigenvalue of $\Delta_2(Q^\pi)$ equals $S_2(G)$.
\end{cor}

The following result which can be found in matrix theory which plays an important role in our proofs. For a symmetric matrix $A$ of order $n$, let $\lambda_1(A)\geq\cdots\geq\lambda_n(A)$ be its eigenvalues.

\begin{lemma}[\cite{KF}]\label{KF}
Let $A$ and $B$ be two real symmetric matrices of size $n$. Then for any $1\leq k\leq n$,
$$\sum^{k}_{i=1}\lambda_{i}(A+B)\leq\sum^{k}_{i=1}\lambda_{i}(A)+\sum^{k}_{i=1}\lambda_{i}(B).$$
\end{lemma}

Apparently, $\lambda_{1}(\Delta_k(A))=\sum^{k}_{i=1}\lambda_{i}(A)$,
then the above lemma can be directly obtained from the fact that
$\lambda_{1}(\Delta_k(A+B))\leq\lambda_{1}(\Delta_k(A))+\lambda_{1}(\Delta_k(B)).$ An immediate result of Lemma \ref{KF} is the following corollary which will be used frequently.

\begin{cor}\label{union}
Let $G$ be a graph of order $n$ such that there exist edge disjoint subgraphs of $G$, say $G_1,\cdots,G_r$, with $E(G)=\bigcup^r_{i=1}E(G_i)$. Then $S_k(G)\leq\sum^{r}_{i=1}S_k(G_i)$ for any integer $k$ with $1\leq k\leq n$.
\end{cor}

The following result is also straightforward.

\begin{lemma}\label{connected}
For some $k$, if Theorem\ref{them1} holds for $G$ and $H$, then it does for $G\cup H$.
\end{lemma}

The next lemma is the key to our approach. It gives a sufficient condition for the result of Theorem\ref{them1}, that holds for almost all graphs.

\begin{lemma}\label{ksubgraph}
If $G$ is a graph with a nonempty subgraph $H$ for which $S_k(H)\leq e(H)$, then $S_k(G)<e(G)$.
\end{lemma}

\begin{proof}
Assume that $G$ is a counterexample with a minimum possible number of edges.
By Corollary \ref{union}, we have $e(G)\leq S_k(G)\leq S_k(H)+S_k(G-H)$.
This implies that $S_k(G-H)\geq e(G-H)$,
since $e(G)\geq e(H)+e(G-H)$ and $S_k(H)\leq e(H)$,
which contradicts the minimality of $e(G)$.
\end{proof}

\begin{cor}\label{subgraph}
If $G$ is a graph with a nonempty subgraph $H$ for which $S_2(H)\leq e(H)$, then $S_2(G)<e(G)$ and $f(G)>3$.
\end{cor}



\begin{lemma}[\cite{Amaro}]\label{v10}
	For any graph $G$ with at most 10 vertices, $f(G)\geq f(K^+_{1,n-1})$ with equality if and only if $G\cong K^+_{1,n-1}$.
\end{lemma}

\begin{lemma}[\cite{DCve}]\label{knsn}
	Let $n$ be a positive integer.
	\begin{itemize}
		\item[(1)] The signless Laplacian eigenvalues of the complete graph  $K_n$ are $2n-2$ and $n-2$ with multiplicities 1 and $n-1$, respectively.
		\item[(2)] The signless Laplacian eigenvalues of the star $K_{1,n-1}$ are $n$, 1 and 0 with multiplicities 1, $n-2$, and 1, respectively.
	\end{itemize}
\end{lemma}

\begin{lemma}[\cite{Fritscher}]\label{Ltree}
Let $T$ be a tree with $n$ vertices and let $1\leq k\leq n$.
Then the sum of the $k$ largest Laplacian eigenvalues of $T$ satisfies
$\sum_{i=1}^k\mu_i(G)\leq (n-1)+2k-1-\frac{2k-2}{n}$.
Moreover, equality is achieved only when $k=1$ and $T$ is a star on $n$ vertices.
\end{lemma}

Using the fact that signless Laplacian matrix and Laplacian matrix are similar for a bipartite graph, from the above lemma, we have the following result.

\begin{cor}\label{tree}
Let $G$ be a tree with $n$ vertices, then $f(G)\geq \frac{2}{n}$.
\end{cor}


\begin{lemma}[\cite{DCve}]\label{interlace}
Let $G$ be a graph with $n$ vertices and let $G^\prime$ be a graph obtained from $G$ by inserting a new edge into $G$. Then the signless Laplacian eigenvalues of $G$ and $G^\prime$ interlace, that is,
$$q_1(G^\prime)\geq q_1(G)\geq\cdots\geq q_n(G^\prime)\geq q_n(G).$$
\end{lemma}




\begin{lemma}\label{uv}
	Let $G$ be a graph with $n$ vertices.
	Suppose that there exist two non-adjacent vertices $u,v\in V(G)$ such
	that $q_k(G)\geq d(u)+d(v)+2$ for some integer $k$, $1\leq k\leq n$.
	If $G^\prime$ is the graph obtained from $G$ by inserting an edge $e=uv$ into $G$, then $S_k(G^\prime)<S_k(G)+1$.
\end{lemma}

\begin{proof}
	For $i=1,\ldots,n$, define $\sigma_i=q_i(G^\prime)-q_i(G)$.
	By Lemma \ref{interlace}, $\sigma_i\geq0$, for any $i$.
	Let $d_1\geq\cdots\geq d_n$ and $d_1^\prime\geq\cdots\geq d_n^\prime$ be vertex degrees of $G$ and $G^\prime$, respectively.
	Recall that for any graph $H$, considering the trace of the matrix $Q^2(H)$, we have
	$$\sum^{|V(H)|}_{i=1}q^2_i(H)=\sum_{v\in V(H)}d^2(v)+2e(H).$$
	Applying this fact, we have
	
	\begin{equation*}
	\begin{aligned}
	\sum^n_{i=1}q^2_i(G^\prime)
	&=\sum^n_{i=1}d_i^{\prime2}+2e(G^\prime)\\
	&=\sum^n_{i=1}d_i^2+2e(G)+2d(u)+2d(v)+4\\
	&=\sum^n_{i=1}q^2_i(G)+2(d(u)+d(v)+2).
	\end{aligned}
	\end{equation*}
	
	This yields that
	
	\begin{equation*}
	\begin{aligned}
	2q_k(G)\sum^k_{i=1}\sigma_i
	&\leq\sum^k_{i=1}2\sigma_iq_i(G)\\
	&<\sum^n_{i=1}q^2_i(G^\prime)-\sum^n_{i=1}q^2_i(G)\\
	&=2(d(u)+d(v)+2).
	\end{aligned}
	\end{equation*}
	
	Since $q_k(G)\geq d(u)+d(v)+2$,
	$S_k(G^\prime)-S_k(G)=\sum^k_{i=1}\sigma_i<1$, and the result follows.
\end{proof}

In particular, if $S_k(G)<e(G)+\binom{k+1}{2}$, then $S_k(G^\prime)<e(G^\prime)+\binom{k+1}{2}$.

\section{Proof of Theorem\ref{them1}}

Amaro et al.\cite{Amaro} pointed out that Theorem \ref{them1} holds for all graphs with at most 10 vertices by computational checking.  Assuming   the graph we are considering has at least     11 vertices in the next. We firstly present Lemma \ref{sn+} which indicates the graph $K^+_{1,n-1}$ is the candidate to be extremal to $f(G)$.

\begin{lemma}\label{sn+}
Let $G$ be isomorphic to $K^+_{1,n-1}$. Then $\frac{1.3}{n}<f(G)<\frac{1.5}{n}$.
\end{lemma}

\begin{proof}
The graph $G$ has an obvious partition $G=\bigvee_H\{G_1,G_2,G_3\}$, and the corresponding quotient matrix is

\begin{align*}
Q^\pi=
\left(
\begin{array}{ccc}
3&1&0\\
2&n-1&n-3\\ 
0&1&1
\end{array}
\right).
\end{align*}

By Lemma \ref{join} and direct calculation, the characteristic polynomial of $Q(G)$ is
\begin{align*} 
P_Q(G,x)=P(Q^\pi,x)(x-1)^{n-3}.
\end{align*}

The 2-th additive compound matrix of $Q^\pi$ is

\begin{align*}
\Delta_2(Q^\pi)=
\left(
\begin{array}{ccc}
n+2&n-3&0\\
1&4&1\\ 
0&2&n
\end{array}
\right),
\end{align*}
and its characteristic polynomia  is $P(\Delta_2(Q^\pi),x)=x^3-(6+2n)x^2+(n^2+9n+9)x-3n^2-9n+4$.
By direct calculation we obtain
\begin{align*}
P(\Delta_2(Q^\pi),n+3-\frac{1.3}{n})&>0,\\
P(\Delta_2(Q^\pi),n+3-\frac{1.5}{n})&<0,\\
P(\Delta_2(Q^\pi),3)&>0,\\
P(\Delta_2(Q^\pi),0)&<0.
\end{align*}

Since $q_2\geq2$, by Lemma \ref{interlace},
$P(Q^\pi,x)$ contains the largest two roots of $P_Q(G,x)$.
From Corollary \ref{c2} and the fact that $e(G)=n$, we have $\frac{1.3}{n}<f(G)<\frac{1.5}{n}$.
\end{proof}

By  Lemma \ref{connected}, we may assume that $G$ is a connected graph with at least 11 vertices.
Noting that for $H=4K_2$ or $H=3K_{1,2}$, one has $S_2(H)=e(H)$,
using Corollary \ref{subgraph}, we may assume $G$ contains neither $H=4K_2$ nor $H=3K_{1,2}$ as a subgraph.
It is sufficient to consider  graphs $G$ whose matching number $m(G)$ is at most 3.
In the following, we prove Theorem \ref{them1} for $m(G)=1,2$ and 3, respectively.

\begin{lemma}\label{m1}
Let $G$ be a graph with $m(G)=1$. Then $f(G)>\frac{1.5}{n}$.
\end{lemma}

\begin{proof}
Let $n=|V(G)|$. Since $m(G)=1$, it is easy to check that either $G= K_{1,k-1}\cup (n-k)K_1$ for some $k$, $1\leq k\leq n$ or $G=K_3\cup (n-3)K_1$.
By Lemma \ref{knsn}, the assertion holds.
\end{proof}

Let  $X, \ Y \subset V(G)$, and 
  $E(X,Y)$ ($e(X,Y)$) be the set of  edges (the number of edges) whose one endpoint is in $X$ and one endpoint is in $Y$.
A pendant neighbor of $G$ is a vertex adjacent to a pendant vertex. Let $p(G)$ be the number of pendant vertices of $G$ and $q(G)$ be the number of pendant neighbors.
\newpage 
\begin{lemma}\label{m2}
Let $G$ be a graph with $m(G)=2$. Then $f(G)\geq f(K^+_{1,n-1})$ with equality if and only if $G\cong K^+_{1,n-1}$.
\end{lemma}

\begin{proof}
We may assume that $G$ is a connected graph with at least 11 vertices.
First suppose that $G$ has no $K_3$ as a subgraph,
if $G$ is a tree, by Corollary \ref{tree}, we have $f(G)>\frac{1.5}{n}$.
Since $m(G)=2$, $G$ is isomorphic to $G_5$ in Fig. \ref{fig1},
using the 2-th additive compound matrix approach, we have $f(G)>0.5>\frac{1.5}{n}$, details are deferred in the Appendix.

Next suppose that $G$ has a subgraph $H=K_3$ with $V(H)=\{u,v,w\}$.

If every edge of $G$ has at least one endpoint in $V(H)$, since $m(G)=2$, we only need to consider $q(G)=0,1$ and 2, respectively.
Then $G$ is isomorphic to $G_1\ (1\leq t_1\leq n-2,\ 0\leq t_2,\ t_3\leq n-3)$ or $G_2\ (7\leq t=n-4)$ in Fig. \ref{fig1}.

If there exists an edge $e=ab$ has no endpoint in $V(H)$,
let $M=V(G)-\{a,b,u,v,w\}$,
since $m(G)=2$, $M$ is an independent set and $e(V(H),M)=0$,
and all vertices in $M$ are adjacent to one of the endpoints of $e$, say $a$.
Hence there are no edges between $b$ and $V(H)$, we only need to consider $|N(a)\cap V(H)|=1$ and 3, respectively (isomorphic to $G_2$ when $|N(a)\cap V(H)|=2$).
Then $G$ is isomorphic to $G_3\ (7\leq t=n-4)$ or $G_4\ (7\leq t=n-4)$ in Fig. \ref{fig1}.

\begin{figure}[htbp!]
\centering
\includegraphics[width=14cm,page={1}]{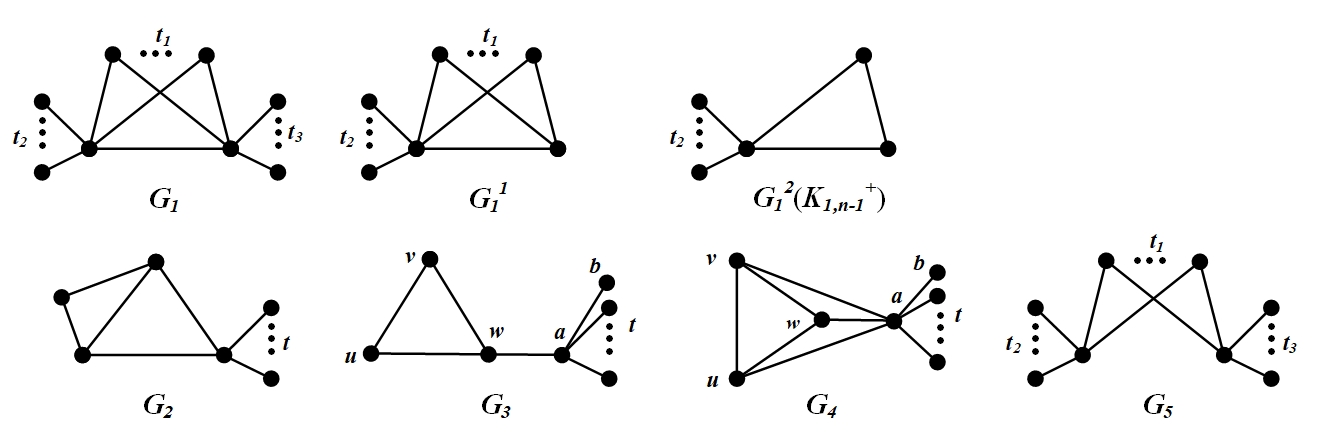}
\caption{Possible forms of graphs $G$ with $m(G)=2$}
\label{fig1}
\end{figure}

When $G\cong G_1$,
let $V(\overline{K_{t_i}})=\{v_1,\ldots,v_{t_i}\}$ $(i=1,2,3)$,
and  $G_{t_1,t_2,t_3}$ be a graph of order $n$ $(n=t_1+t_2+t_3)$,
which consists of three independent sets $\overline{K_{t_i}}$ $(i=1,2,3)$,
where $t_1$ triangles share one common edge and $t_2,t_3$ pendant vertices at each end of the common edge.
The graph $G_{t_1,t_2,t_3}$ has an obvious partition $G=\bigvee_H\{G_1,G_2,...,G_k\}(k=5)$ where the corresponding quotient matrix is

\begin{align*}
Q^\pi=
\left(
\begin{array}{ccccc}
2&0&0&1&1\\
0&1&0&1&0\\ 
0&0&1&0&1\\
t_1&t_2&0&t_1+t_2+1&1\\
t_1&0&t_3&1&t_1+t_3+1
\end{array}
\right).
\end{align*}

By Lemma \ref{join} and direct calculation, the characteristic polynomial of $Q(G_{t_1,t_2,t_3})$ is
\begin{align*} 
P_Q(G_{t_1,t_2,t_3},x)=P(Q^\pi(G_{t_1,t_2,t_3}),x)(x-2)^{t_1-1}(x-1)^{t_2+t_3-2}.
\end{align*}

Without loss of generality, we may assume   $t_2\geq t_3$.
 Let
\begin{align*} 
P_1&=P(\Delta_2(Q^\pi(G_{t_1,t_2,t_3})),x+e(G_{t_1,t_2,t_3})+3),\\
P_2&=P(\Delta_2(Q^\pi(G_{t_1,t_2+1,t_3-1})),x+e(G_{t_1,t_2+1,t_3-1})+3).
\end{align*}

By direct computation, $P_1-P_2$ is always positive
(details are deferred in the Appendix),
which indicates the largest root of $P_1$ less than the largest root of $P_2$.

Since $q_2>2$, by Lemma \ref{interlace}
(if $t_1=1$ and $t_3=0$, the previous results (see Theorem 2.4 in \cite{Aouchiche}) verify that $q_2>3-\frac{2.5}{n}>2$),
$P(Q^\pi(G_{t_1,t_2,t_3},x))$ and $P(Q^\pi(G_{t_1,t_2+1,t_3-1},x))$ contain  the largest two roots of $P_Q(G_{t_1,t_2,t_3},x)$ and $P_Q(G_{t_1,t_2+1,t_3-1},x)$.
By Corollary \ref{c2}, we have $f(G_{t_1,t_2,t_3})>f(G_{t_1,t_2+1,t_3-1})$.
By recursive relationships, we have $f(G_1)>f(G_1^1)$.

Let $G_{t_1,t_2}$ be a graph of order $n$ $(n=t_1+t_2)$,
which consists of two independent sets $\overline{K_{t_i}}$ $(i=1,2)$,
where $t_1$ triangles share one common edge and $t_2$ pendant vertices at the end of the common edge.
Let
\begin{align*} 
P_3&=P(\Delta_2(Q^\pi(G_{t_1,t_2})),x+e(G_{t_1,t_2})+3),\\
P_4&=P(\Delta_2(Q^\pi(G_{t_1-1,t_2+1})),x+e(G_{t_1-1,t_2+1})+3).
\end{align*}

By using computer computations, $P_3-P_4$ is always positive
(details are deferred in the Appendix),
which indicates the largest root of $P_3$ less than the largest root of $P_4$.

Since $q_2>2$, by Lemma \ref{interlace}
(if $t_1=1$, the previous results (see Theorem 2.4 in \cite{Aouchiche}) verify that $q_2>3-\frac{2.5}{n}>2$),
$P(Q^\pi(G_{t_1,t_2},x))$ and $P(Q^\pi(G_{t_1-1,t_2+1},x))$ contains the largest two roots of $P_Q(G_{t_1,t_2},x)$ and $P_Q(G_{t_1-1,t_2+1},x)$.
By Corollary \ref{c2}, we have $f(G_{t_1,t_2})>f(G_{t_1-1,t_2+1})$.
By recursive relationships, we have $f(G_1^1)>f(G_1^2)$.
Then $f(G_1)>f(G_1^1)>f(G_1^2)=f(K^+_{1,n-1})$.

When $G\cong G_2$,
the graph $G$ has an obvious partition $G=\bigvee_H\{G_1,G_2,...,G_k\}(k=4)$ where the corresponding quotient matrix is

\begin{align*}
Q^\pi=
\left(
\begin{array}{cccc}
1&1&0&0\\ 
t&t+2&2&0\\
0&1&4&1\\
0&0&2&2
\end{array}
\right).
\end{align*}

By Lemma \ref{join} and direct calculation, the characteristic polynomial of $Q(G)$ is
\begin{align*} 
P_Q(G,x)=P(Q^\pi,x)(x-2)(x-1)^{t-1}.
\end{align*}

By simple calculations we obtain
\begin{align*} 
P(Q^\pi,t+3.25-\frac{1.5}{t+4})&=\frac{0.25}{(4+t)^4}g(t)>0,\\
P(Q^\pi,4.75)&=-0.3t-19.98<0,\\
P(Q^\pi,2)&=4t>0,\\
P(Q^\pi,1)&=-t<0,\\
P(Q^\pi,0)&=8>0,
\end{align*}
with $g(t)=(t^2+5.3t+8.9)(t+5.2)(t+4.6)(t+4.3)(t+1.3)(t-6.9)$,
which indicates $q_1(G)<t+3.25-\frac{1.5}{t+4}$ and $q_2(G)<4.75$, then $S_2(G)<t+8-\frac{1.5}{t+4}=e(G)+3-\frac{1.5}{n}$,  we have $f(G)>\frac{1.5}{n}$.

When $G\cong G_3$, we use the quotient matrix approach,
when $G\cong G_4$, we use the 2-th additive compound matrix approach,
details are deferred in the Appendix.
This completes the proof.

\end{proof}

For nonempty set  $X  \subset V(G)$, let 
$N(X)=\{ v\mid v\not\in X \ \text{and} \  N(v)\cap X\neq \emptyset\}$. Suppose that $V^\prime$, $E^\prime$ are nonempty subsets of $V(G)$ and $E(G)$, the subgraphs induced by  $V^\prime$ and  $E^\prime$ are  denoted by $G[V^\prime]$ and  $G[E^\prime]$, respectively. 
\begin{lemma}\label{m3}
Let $G$ be a graph with $m(G)=3$. Then $f(G)>\frac{1.5}{n}$.
\end{lemma}

\begin{proof}
Without loss of generality, we may suppose that $e_1=a_1b_1$, $e_2=a_2b_2$ and $e_3=a_3b_3$ are three independent edges in $G$.

We first suppose that $G$ contains $K_3\cup2K_2$ as a subgraph and $V(K_3)=\{a_1,b_1,u\}$.
Let $V_1=V(K_3\cup2K_2)$, $M^\prime=V(G)-V_1$. Then $|M^\prime|\geq4$.
Since $m(G)=3$, $M^\prime$ is an independent set and $e(V(K_3),M^\prime)=0$.
Since $G$ has neither $3K_{1,2}$ nor $4K_2$ as a subgraph and $|M^\prime|\geq4$,
there exists exactly one vertex $a_3$ (w.l.o.g.) which is adjacent to vertices of $M^\prime$ and $d(b_3)=1$.

Let $V_2=\{a_1,b_1,u,a_2,b_2\}$, $E_1=\{a_3w|w\in V_2\}$.
Then $G[E_1]$ is a star with center $a_3$, and $G-E_1$ is a disjoint union of a star $K_{1,|M^\prime|+1}$ with center $a_3$ and a graph $G[V_2]$ containing $K_3+K_2$.

If $G[V_2]$ is connected, then $G[V_2]$ is spanning subgraph of $K_5$ which contains $K_3\cup K_2$ as a subgraph, by direct calculation, we have $q_1(G[V_2])<e(G[V_2])-0.15<e(G[V_2])-\frac{1.5}{n}\ (n\geq11)$.
Since $n(K_{1,|M^\prime|+1})\geq6$, so $q_2(G[V_2])\leq q_2(K_5)<q_1(K_{1,|M^\prime|+1})$.
Hence $S_2(G-E_1)=q_1(G[V_2])+q_1(K_{1,|M^\prime|+1})<e(G[V_2])-\frac{1.5}{n}+e(K_{1,|M^\prime|+1})+1=e(G-E_1)+1-\frac{1.5}{n}$.
Thus by Corollary \ref{union}, $S_2(G)\leq S_2(G[E_1])+S_2(G-E_1)<e(G[E_1])+2+e(G-E_1)+1-\frac{1.5}{n}=e(G)+3-\frac{1.5}{n}$,
we have $f(G)>\frac{1.5}{n}$.

If $G[V_2]$ is disconnected, then $G[V_2]=K_3\cup K_2$.
Let $t=d(a_2)+d(b_2)$. We have $t=3, 4$.
It is not hard to see that $m(G-e_2)=2$ and $G-e_2$ contains disjoint union of a star $K_{1,|M^\prime|+1}$ with center $a_3$ and $K_3$.
Therefore, by Lemma \ref{interlace}, $q_2(G-e_2)\geq t$.
Using Lemmas \ref{m2} and  \ref{uv}, we find that $S_2(G)<S_2(G-e_2)+1<(e(G-e_2)+3-\frac{1.5}{n})+1=e(G)+3-\frac{1.5}{n}$,
we have $f(G)>\frac{1.5}{n}$.

Next suppose that $G$ has no subgraph $K_3\cup2K_2$.
Since $m(G)=3$, $M=V(G)-V(\{e_1,e_2,e_3\})$ is an independent set and $|M|\geq5$.
Since $G$ has no $4K_2$ and $K_3+2K_2$ as subgraph,
either $N_M(a_i)=\emptyset$ or $N_M(b_i)=\emptyset$ for $i=1,2,3$.
Without loss of generality, we may assume that $N(M)\subseteq \{a_1,a_2,a_3\}$.
We consider the following three cases.

\textbf{Case 1:} $|N(M)|=3$.\par
We have $N(M)=\{a_1,a_2,a_3\}$.
The bipartite subgraph $G[E(\{a_1,a_2,a_3\},M)]$ has no perfect matching since $G$ has no $3K_{1,2}$.
By Hall’s Theorem, there exists a subset of $\{a_1,a_2,a_3\}$ with 2 elements,
say $\{a_2,a_3\}$, such that $|N(\{a_2,a_3\})\cap M|=1$
(we have $|N(\{a_2,a_3\})\cap M|\neq0$ since $|N(M)|=3$).
This means that there exists exactly one vertex $y\in M$ which is adjacent to both $a_2$ and $a_3$.
If $d(b_1)\geq2$, then we clearly find a subgraph isomorphic to $3K_{1,2}$ in $G$, a contradiction.
Therefore, $d(b_1)=1$.
Let $V_1=\{a_2,b_2,a_3,b_3,y\}$, $E_1=\{a_1w|w\in V_1\}$.
Then $G[E_1]$ is a star with center $a_1$, and $G-E_1$ is a disjoint union of a star $K_{1,|M|}$ with center $a_1$ and a graph $G[V_1]$ containing $P_5$.
Since G has no $K_3\cup2K_2$, two edges $a_2a_3$ and $b_2b_3$ can not exist at the same time,
therefore $G[V_1]$ is isomorphic to graphs in Fig. \ref{fig2}.

\begin{figure}[htbp!]
\centering
\includegraphics[width=8cm,page={1}]{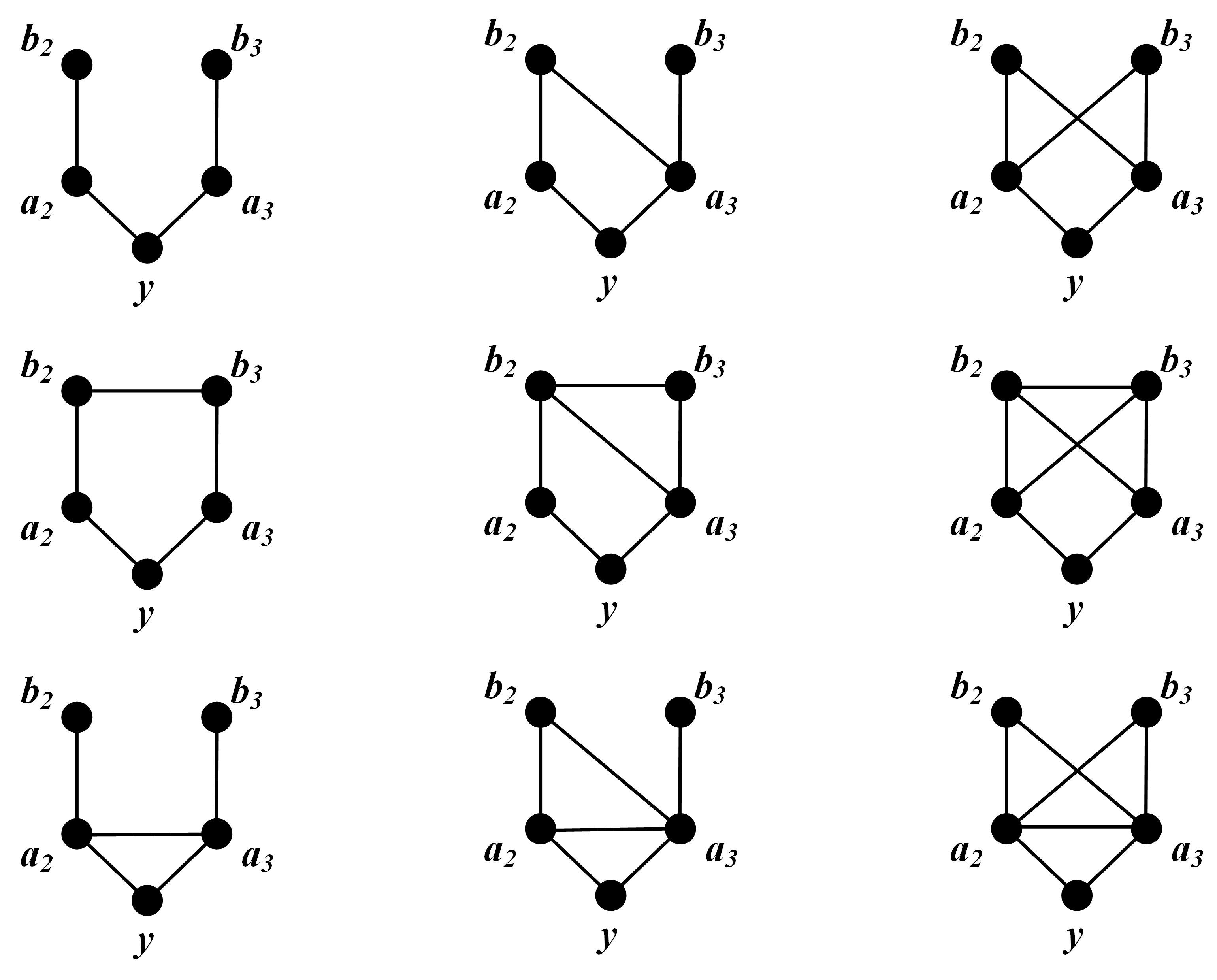}
\caption{Possible forms of graphs containing a $P_5$}
\label{fig2}
\end{figure}

By direct calculation, we have $q_1(G[V_1])<e(G[V_1])-\frac{1.5}{n}\ (n\geq11)$.
Since $n(K_{1,|M|})\geq6$, so $q_2(G[V_1])\leq q_1(K_{1,|M|})$.
Hence $S_2(G-E_1)=q_1(G[V_1])+q_1(K_{1,|M|})<e(G[V_1])-\frac{1.5}{n}+e(K_{1,|M|})+1=e(G-E_1)+1-\frac{1.5}{n}$.
Thus by Corollary \ref{union}, $S_2(G)\leq S_2(G[E_1])+S_2(G-E_1)<e(G[E_1])+2+e(G-E_1)+1-\frac{1.5}{n}=e(G)+3-\frac{1.5}{n}$,
we have $f(G)>\frac{1.5}{n}$.

\textbf{Case 2:} $|N(M)|=2$.\par
Without loss of generality, suppose that $N(M)=\{a_1,a_2\}$.
Since $m(G)=3$, $b_1$ is not adjacent to $b_2$.
If $b_1$ is adjacent to $a_3$ or $b_3$, then changing the role of $e_1$, $e_2$, $e_3$ by three independent edges $\{a_1,z\}$, $e_2$, $e_3$ for some vertex $z\in M\cap N(a_1)$, we have Case 1.
Therefore, we may assume that $b_1$, and similarly $b_2$, is adjacent to none of the vertices $a_3$ and $b_3$.
Hence we only need to consider $|N(\{a_3,b_3\})|=1$ or 2, respectively. Consider the existence of the edge $a_1a_2$, possible forms of graph $G$ are in Fig. \ref{fig3}.

\begin{figure}[htbp!]
\centering
\includegraphics[width=14cm,page={1}]{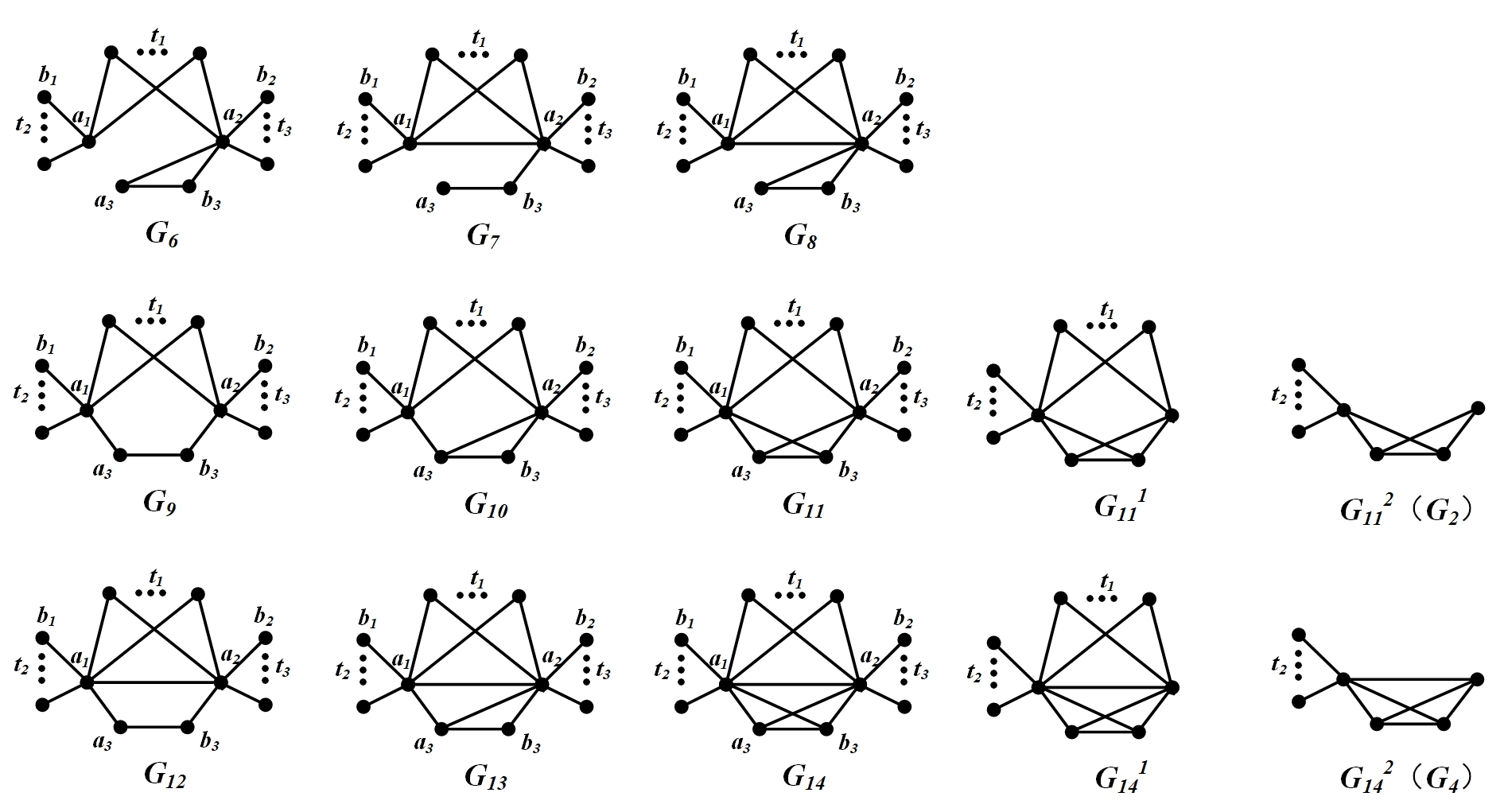}
\caption{Possible forms of graphs $G$ with $m(G)=3$ and $|N(M)|=2$}
\label{fig3}
\end{figure}

When $G\cong G_i$, $i=6,7,8,9,10,12,13$,
let $t=d(a_3)+d(b_3)$, we have $t\in \{3,4,5\}$.

We first present several special cases:
$t_1=t_2=1$ when $G\cong G_6$; $t_1=0$ and $t_2=2$; $t_1=t_2=1$; $t_1=2$ and $t_2=0$ when $G\cong G_8,G_{10},G_{13}$
(were checked by computer computations).

Except these special cases,
it is not hard to see that $m(G-e_3)=2$ and $G-e_3$ contains two disjoint stars $K_{1,t-1}$ with centers $a_1$ and $a_2$
Therefore, by Lemma \ref{interlace}, $q_2(G-e_3)\geq t$.
Using Lemmas \ref{m2} and  \ref{uv}, we find that $S_2(G)<S_2(G-e_3)+1<(e(G-e_3)+3-\frac{1.5}{n})+1=e(G)+3-\frac{1.5}{n}$,
we have $f(G)>\frac{1.5}{n}$.

When $G\cong G_{11}$, we use the 2-th additive compound matrix approach,
by using computer computations and recursive relationships,
we have $f(G_{11})>f(G_{11}^1)$ and $f(G_{11}^1)>f(G_{11}^2)$.
Then $f(G_{11})>f(G_{11}^1)>f(G_{11}^2)=f(G_2)$,
we have $f(G)>\frac{1.5}{n}$.

When $G\cong G_{14}$, we use the 2-th additive compound matrix approach,
by using computer computations and recursive relationships,
we have $f(G_{14})>f(G_{14}^1)$ and $f(G_{14}^1)>f(G_{14}^2)$.
Then $f(G_{14})>f(G_{14}^1)>f(G_{14}^2)=f(G_4)$,
we have $f(G)>\frac{1.5}{n}$.

\textbf{Case 3:} $|N(M)|=1$.\par
Without loss of generality, suppose that $N(M)=\{a_1\}$.
If $d(b_1)\geq2$, then we clearly find three independent edges $e^\prime_1$, $e^\prime_2$, $e^\prime_3$ in $G$ such that the set $M^\prime=V(G)-V(\{e^\prime_1, e^\prime_2, e^\prime_3\}$ is an independent set and $|N(M^\prime)|\geq2$ which is dealt with as the previous cases.
Therefore, $d(b_1)=1$.
Let $V_1=\{a_2,b_2,a_3,b_3\}$, $E_1=\{a_1w|w\in V_1\}$.
Then $G[E_1]$ is a star with center $a_1$, and $G-E_1$ is a disjoint union of a star $K_{1,|M|+1}$ with center $a_1$ and a graph $G[V_1]$ containing $2K_2$.
Therefore $G[V_1]$ is isomorphic to graphs in Fig. \ref{fig4}.

\begin{figure}[htbp!]
\centering
\includegraphics[width=8cm,page={1}]{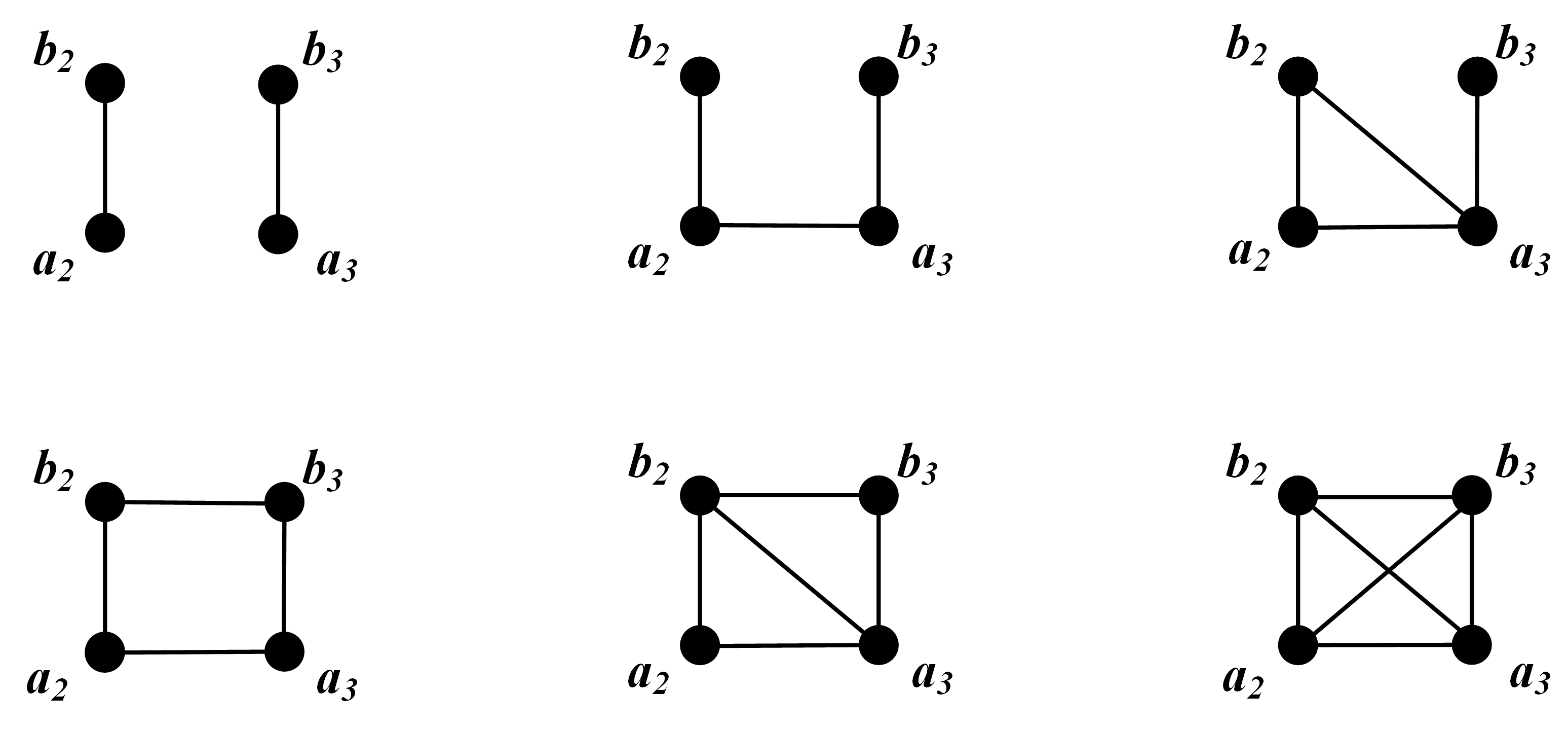}
\caption{Possible forms of graphs containing $2K_2$}
\label{fig4}
\end{figure}

By direct calculation, we have $q_1(G[V_1])<e(G[V_1])+1-\frac{1.5}{n}\ (n\geq11)$,
hence $S_2(G-E_1)=q_1(G[V_1])+q_1(K_{1,|M|+1})<e(G[V_1])+1-\frac{1.5}{n}+e(K_{1,|M|})+1=e(G-E_1)+2-\frac{1.5}{n}$.
If $|N(a_1)\cap V_1|=1$, which means $|E_1|=1$, by Corollary \ref{union}, $S_2(G)\leq S_2(G[E_1])+S_2(G-E_1)<e(G[E_1])+1+e(G-E_1)+2-\frac{1.5}{n}=e(G)+3-\frac{1.5}{n}$,
we have $f(G)>\frac{1.5}{n}$.
Hence we only need to consider $|N(a_1)|\cap V_1|=2$, 3 and 4, respectively.
Possible forms of graph $G$ are in Fig. \ref{fig5} (If $G$ is a tree graph, the assertion holds by Corollary \ref{tree}).

\begin{figure}[htbp!]
\centering
\includegraphics[width=13cm,page={1}]{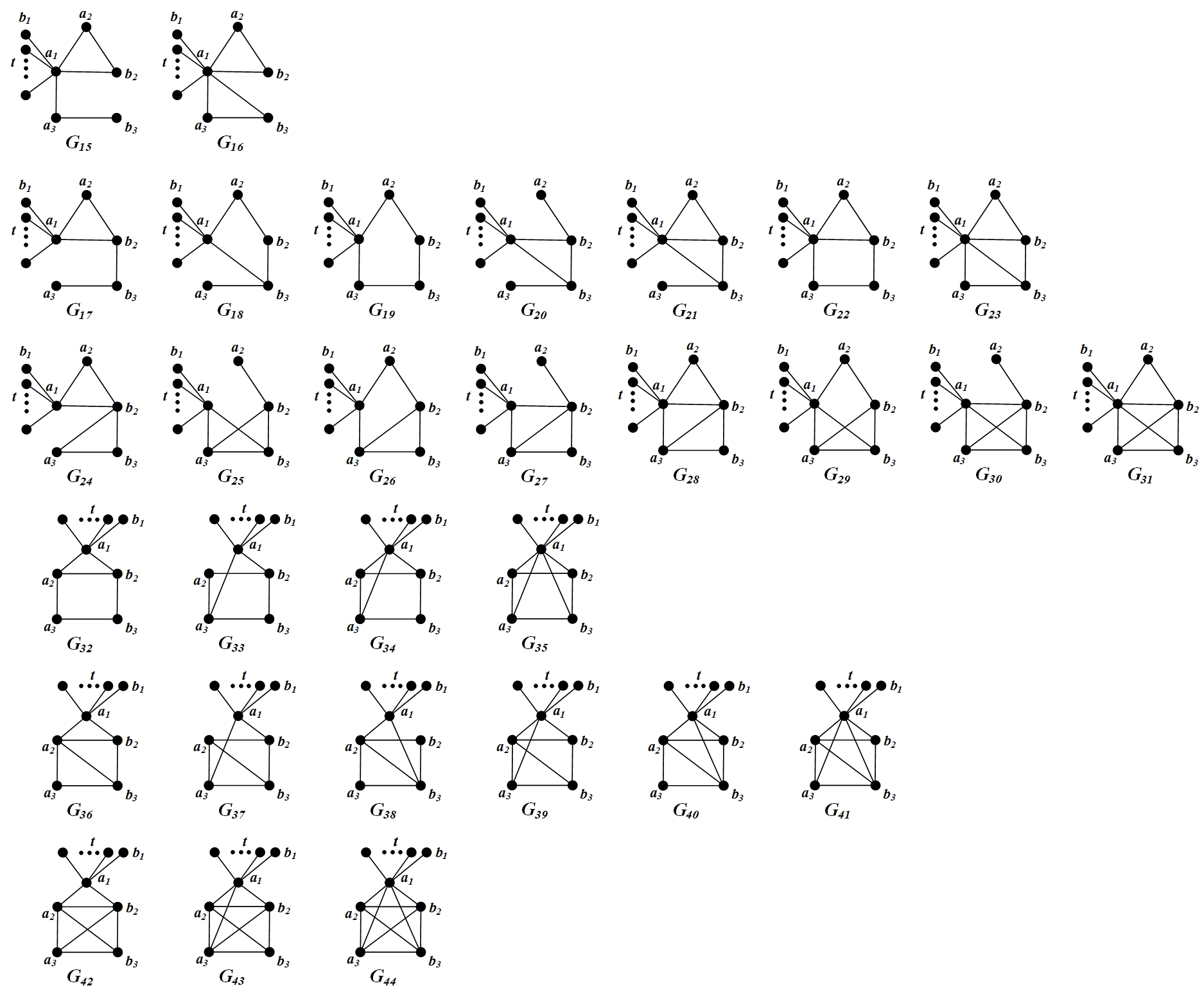}
\caption{Possible forms of graphs $G$ with $m(G)=3$ and $|N(M)|=1$}
\label{fig5}
\end{figure}

When $G\cong G_i$, $i=17,24,25$,
let $t=d(a_3)+d(b_3)$ (pick $a_2$ and $b_2$ when $G\cong G_{25}$), we have $t\in \{3,4 \}$.
It is not hard to see that $m(G-e_3)=2$ and $G-e_3$ contains two disjoint stars $K_{1,t-1}$ with centers $a_1$ and $b_2$ (contains disjoint union of a star $K_{1,t-1}$ with center $a_1$ and a $K_3$ when $G\cong G_{25}$).
Therefore, by Lemma \ref{interlace}, $q_2(G-e_3)\geq t$.
Using Lemmas \ref{m2} and  \ref{uv}, we find that $S_2(G)<S_2(G-e_3)+1<(e(G-e_3)+3-\frac{1.5}{n})+1=e(G)+3-\frac{1.5}{n}$,
we have $f(G)>\frac{1.5}{n}$.

When $G$ is isomorphic to other graphs in Fig. \ref{fig5}, we use the quotient matrix approach, details are deferred in the Appendix.
This completes the proof.
\end{proof}

\noindent\textbf{Proof of Theorem \ref{them1}:} The result follows immediately from Corollary \ref{subgraph}, Lemmas \ref{sn+}, \ref{m1}, \ref{m2} and \ref{m3}.

\section{Declaration of competing interest}

There is no competing interest.

\newpage
\section*{Appendix}
The details of polynomial $P_1-P_2$ and $P_3-P_4$ can be obtained from:

https://github.com/ZZM0329/Polynomial/tree/main/G1

    \centering
    \begin{longtable}{ccl}
    \caption{2-th additive compound matrix approach}\\
    \hline
No. &Corresponding quotient matrix $Q^\pi$ &Characteristic polynomial $P_Q(G,x)$\\
    \hline
        $G_4$ & 
$\left(\begin{array}{ccc}
1&1&0\\ 
t&t+3&3\\
0&1&5
\end{array}\right)$ &       
         $P(Q^\pi,x)(x-2)^2(x-1)^{t-1}$ \\
    \hline
         $G_5$ & 
$\left(\begin{array}{ccccc}
2&0&0&1&1\\
0&1&0&1&0\\
0&0&1&0&1\\
t_1&t_2&0&t_1+t_2&0\\
t_1&0&t_3&0&t_1+t_3
\end{array}\right)$ &       
         $P(Q^\pi,x)(x-2)^{t_1-1}(x-1)^{t_2+t_3-2}$ \\
    \hline
         $G_{11}$ & 
$\left(\begin{array}{cccccc}
2&0&0&1&1&0\\
0&1&0&1&0&0\\
0&0&1&0&1&0\\
t_1&t_2&0&t_1+t_2+2&0&2\\
t_1&0&t_3&0&t_1+t_3+2&2\\
0&0&0&1&1&4
\end{array}\right)$ &       
         $P(Q^\pi,x)(x-2)^{t_1}(x-1)^{t_2+t_3-2}$ \\
    \hline        
         $G_{14}$ & 
$\left(\begin{array}{cccccc}
2&0&0&1&1&0\\
0&1&0&1&0&0\\
0&0&1&0&1&0\\
t_1&t_2&0&t_1+t_2+3&1&2\\
t_1&0&t_3&1&t_1+t_3+3&2\\
0&0&0&1&1&4
\end{array}\right)$ &       
         $P(Q^\pi,x)(x-2)^{t_1}(x-1)^{t_2+t_3-2}$ \\
    \hline
    \end{longtable}

    \centering
    \begin{longtable}{ccp{7.5cm}}
    \caption{Quotient matrix approach}\\
    \hline
         No. & Corresponding quotient matrix $Q^\pi$ & Characteristic polynomial $P_Q(G,x)$\\
    \hline
         $G_3$ & 
$\left(\begin{array}{cccc}
1&1&0&0\\ 
t&t+1&1&0\\
0&1&3&2\\
0&0&1&3
\end{array}\right)$ &       
         $[x^4-(t+8)x^3+(6t+19)x^2-(7t+16)x+4](x-1)^t$ \\
    \hline    
         $G_{15}$ & 
$\left(\begin{array}{ccccc}
1&1&0&0&0\\ 
t&t+3&2&0&1\\
0&1&3&0&0\\
0&0&0&1&1\\
0&1&0&1&2
\end{array}\right)$ &       
         $[x^5-(t+10)+x^4+(6t+34)x^3-(10t+48)x^2+(3t+27)x-4](x-1)^t$ \\
    \hline    
         $G_{16}$ & 
$\left(\begin{array}{cccc}
1&1&0&0\\ 
t&t+4&2&2\\
0&1&3&0\\
0&1&0&3
\end{array}\right)$ &       
         $[x^3-(t+8)x^2+(3t+15)x-8](x-3)(x-1)^{t+1}$ \\
    \hline
         $G_{18}$ & 
$\left(\begin{array}{cccccc}
1&1&0&0&0&0\\
t&t+2&1&0&1&0\\
0&1&2&1&0&0\\
0&0&1&2&1&0\\
0&1&0&1&3&1\\
0&0&0&0&1&1
\end{array}\right)$ &       
         $[x^4-(t+10)x^3+(7t+34)x^2-(13t+46)x+4t+20](x-1)^{t}x$ \\
    \hline        
         $G_{19}$ & 
$\left(\begin{array}{cccccc}
1&1&0&0&0&0\\
t&t+2&1&0&0&1\\
0&1&2&1&0&0\\
0&0&1&2&1&0\\
0&0&0&1&2&1\\
0&1&0&0&1&2
\end{array}\right)$ &       
         $[x^4-(t+8)x^3+(5t+20)x^2-(5t+17)x+4](x^2-3x+1)(x-1)^{t-1}$ \\
    \hline        
         $G_{20}$ & 
$\left(\begin{array}{cccccc}
1&1&0&0&0&0\\
t&t+2&0&1&1&0\\
0&0&1&1&0&0\\
0&1&1&3&1&0\\
0&1&0&1&3&1\\
0&0&0&0&1&1
\end{array}\right)$ &       
         $[x^4-(t+8)x^3+(5t+18)x^2-(3t+15)x+4](x^2-3x+1)(x-1)^{t-1}$ \\
    \hline        
         $G_{21}$ & 
$\left(\begin{array}{cccccc}
1&1&0&0&0&0\\
t&t+3&1&1&1&0\\
0&1&2&1&0&0\\
0&1&1&3&1&0\\
0&1&0&1&3&1\\
0&0&0&0&1&1
\end{array}\right)$ &       
         $[x^6-(t+13)x^5+(9t+62)x^4-(26t+140)x^3+(27t+158)x^2-(8t+84)x+16](x-1)^{t-1}$ \\
    \hline       
         $G_{22}$ & 
$\left(\begin{array}{cccccc}
1&1&0&0&0&0\\
t&t+3&1&1&0&1\\
0&1&2&1&0&0\\
0&1&1&3&1&0\\
0&0&0&1&2&1\\
0&1&0&0&1&2
\end{array}\right)$ &       
         $[x^6-(t+13)x^5+(9t+63)x^4-(27t+145)x^3+(31t+165)x^2-(11t+87)x+16](x-1)^{t-1}$ \\
    \hline       
         $G_{23}$ & 
$\left(\begin{array}{cccccc}
1&1&0&0&0&0\\
t&t+4&1&1&1&1\\
0&1&2&1&0&0\\
0&1&1&3&1&0\\
0&1&0&1&3&1\\
0&1&0&0&1&2
\end{array}\right)$ &       
         $[x^4-(t+11)x^3+(6t+37)x^2-(7t+47)x+20](x-3)(x-1)^t$ \\
    \hline       
         $G_{26}$ & 
$\left(\begin{array}{cccccc}
1&1&0&0&0&0\\
t&t+2&1&0&0&1\\
0&1&2&1&0&0\\
0&0&1&3&1&1\\
0&0&0&1&2&1\\
0&1&0&1&1&3
\end{array}\right)$ &       
         $[x^5-(t+12)x^4+(9t+51)x^3-(24t+94)x^2+(19t+71)x-16](x-1)^t$ \\
    \hline       
         $G_{27}$ & 
$\left(\begin{array}{cccccc}
1&1&0&0&0&0\\
t&t+2&0&1&0&1\\
0&0&1&1&0&0\\
0&1&1&4&1&1\\
0&0&0&1&2&1\\
0&1&0&1&1&3
\end{array}\right)$ &       
         $[x^6-(t+13)x^5+(10t+61)x^4-(31t+135)x^3+(35t+150)x^2-(12t+80)x+16](x-1)^{t-1}$ \\
    \hline       
         $G_{28}$ & 
$\left(\begin{array}{cccccc}
1&1&0&0&0&0\\
t&t+3&1&1&0&1\\
0&1&2&1&0&0\\
0&1&1&4&1&1\\
0&0&0&1&2&1\\
0&1&0&1&1&3
\end{array}\right)$ &       
         $[x^6-(t+15)x^5+(11t+84)x^4-(40t+228)x^3+(58t+319)x^2-(29t+221)x+60](x-1)^{t-1}$ \\
    \hline       
         $G_{29}$ & 
$\left(\begin{array}{ccccc}
1&1&0&0&0\\
t&t+3&1&0&2\\
0&1&2&1&0\\
0&0&1&3&2\\
0&1&0&1&4
\end{array}\right)$ &       
         $[x^5-(t+13)x^4+(9t+59)x^3-(23t+115)x^2+(16t+92)x-24](x-2)(x-1)^{t-1}$ \\
    \hline       
         $G_{30}$ & 
$\left(\begin{array}{ccccc}
1&1&0&0&0\\
t&t+3&0&1&2\\
0&0&1&1&0\\
0&1&1&4&2\\
0&1&0&1&4
\end{array}\right)$ &       
         $[x^5-(t+13)x^4+(9t+57)x^3-(21t+107)x^2+(10t+86)x-24](x-2)(x-1)^{t-1}$ \\
    \hline       
         $G_{31}$ & 
$\left(\begin{array}{ccccc}
1&1&0&0&0\\
t&t+4&1&1&2\\
0&1&2&1&0\\
0&1&1&4&2\\
0&1&0&1&4
\end{array}\right)$ &       
         $[x^4-(t+12)x^3+(7t+43)x^2-(8t+56)x+24](x-3)(x-2)(x-1)^{t-1}$ \\
    \hline           
         $G_{32}$ & 
$\left(\begin{array}{cccccc}
1&1&0&0&0&0\\
t&t+2&1&1&0&0\\
0&1&3&1&0&1\\
0&1&1&3&1&0\\
0&0&0&1&2&1\\
0&0&1&0&1&2
\end{array}\right)$ &       
         $[x^4-(t+10)x^3+(7t+32)x^2-(11t+39)x+16](x^2-3x+1)(x-1)^{t-1}$ \\
    \hline
         $G_{33}$ & 
$\left(\begin{array}{cccc}
1&1&0&0\\ 
t&t+2&0&2\\
0&0&2&2\\
0&1&2&3
\end{array}\right)$ &       
         $[x^3-(t+8)x^2+(5t+17)x-2t-10](x-3)(x-2)(x-1)^{t-1}x$ \\
    \hline
         $G_{34}$ & 
$\left(\begin{array}{ccccc}
1&1&0&0&0\\ 
t&t+3&1&2&0\\
0&1&3&2&0\\
0&1&1&3&1\\
0&0&0&2&2
\end{array}\right)$ &       
         $[x^5-(t+12)x^4+(8t+49)x^3-(17t+86)x^2+(8t+64)x-16](x-3)(x-1)^{t-1}$ \\
    \hline        
         $G_{35}$ & 
$\left(\begin{array}{ccc}
1&1&0\\ 
t&t+4&4\\
0&1&5
\end{array}\right)$ &       
         $[x^3-(t+10)x^2+(5t+25)x-16](x-3)^2(x-1)^t$ \\
    \hline
         $G_{36}$ & 
$\left(\begin{array}{cccccc}
1&1&0&0&0&0\\
t&t+2&1&1&0&0\\
0&1&4&1&1&1\\
0&1&1&3&1&0\\
0&0&1&1&3&1\\
0&0&1&0&1&2
\end{array}\right)$ &       
         $[x^6-(t+15)x^5+(12t+84)x^4-(48t+228)x^3+(77t+319)x^2-(41t+221)x+60](x-1)^{t-1}$ \\
    \hline       
         $G_{37}$ & 
$\left(\begin{array}{cccc}
1&1&0&0\\ 
t&t+2&0&2\\
0&0&4&2\\
0&1&2&3
\end{array}\right)$ &       
         $[x^4-(t+10)x^3+(7t+29)x^2-(8t+28)x+8](x-3)(x-2)(x-1)^{t-1}$ \\
    \hline
         $G_{38}$ & 
$\left(\begin{array}{cccc}
1&1&0&0\\ 
t&t+2&2&0\\
0&1&5&2\\
0&0&2&2
\end{array}\right)$ &       
         $[x^3-(t+9)x^2+(6t+18)x-8](x-3)(x-2)(x-1)^t$ \\
    \hline
         $G_{39}$ & 
$\left(\begin{array}{ccccc}
1&1&0&0&0\\ 
t&t+3&1&2&0\\
0&1&4&2&1\\
0&1&1&3&1\\
0&0&1&2&3
\end{array}\right)$ &       
         $[x^5-(t+14)x^4+(10t+68)x^3-(28t+146)x^2+(23t+139)x-48](x-3)(x-1)^{t-1}$ \\
    \hline        
         $G_{40}$ & 
$\left(\begin{array}{ccccc}
1&1&0&0&0\\ 
t&t+3&2&1&0\\
0&1&5&1&1\\
0&1&2&3&0\\
0&0&2&0&2
\end{array}\right)$ &       
         $[x^5-(t+14)x^4+(10t+67)x^3-(27t+142)x^2+(20t+136)x-48](x-3)(x-1)^{t-1}$ \\
    \hline
         $G_{41}$ & 
$\left(\begin{array}{cccc}
1&1&0&0\\ 
t&t+4&2&2\\
0&1&5&2\\
0&1&2&3
\end{array}\right)$ &       
         $[x^4-(t+13)x^3+(8t+51)x^2-(11t+75)x+36](x-3)^2(x-1)^{t-1}$ \\
    \hline
         $G_{42}$ & 
$\left(\begin{array}{cccc}
1&1&0&0\\ 
t&t+2&2&0\\
0&1&5&2\\
0&0&2&4
\end{array}\right)$ &       
         $[x^4-(t+12)x^3+(9t+43)x^2-(16t+56)x+24](x-3)(x-2)(x-1)^{t-1}$ \\
    \hline
         $G_{43}$ & 
$\left(\begin{array}{ccccc}
1&1&0&0&0\\ 
t&t+3&2&1&0\\
0&1&5&1&1\\
0&1&2&4&1\\
0&0&2&1&3
\end{array}\right)$ &       
         $[x^4-(t+13)x^3+(9t+51)x^2-(15t+75)x+36](x-3)^2(x-1)^{t-1}$ \\
    \hline
         $G_{44}$ & 
$\left(\begin{array}{ccc}
1&1&0\\ 
t&t+4&4\\
0&1&7
\end{array}\right)$ &       
         $[x^3-(t+12)x^2+(7t+35)x-24](x-3)^3(x-1)^{t-1}$ \\    
    \hline
    \end{longtable}


\begin{thebibliography}{3}
	\bibitem{Amaro}B. Amaro, L.S. Lima, C.S. Oliveira, C. Lavor, N. Abreu, A note on the sum of the largest signless Laplacian eigenvalues, \emph{Electronic Notes in Discrete Mathematics.}, \textbf{54}(2016), 175--180.
	
	\bibitem{Aouchiche}M. Aouchiche, P. Hansen, C. Lucas, On the extremal values of the second largest Q-
	eigenvalue, \emph{Linear Algebra Appl.}, \textbf{435}(2011), 2591--2606.
	
	\bibitem{Ashraf}F. Ashraf, G.R. Omidi, B. Tayfeh-Rezaibe, On the sum of signless Laplacian eigenvalues of graph, \emph{Linear Algebra Appl.}, \textbf{438}(2013), 4539--4546.
	
	\bibitem{Brouwer}A.E. Brouwer, W.H. Haemers, Spectra of graphs, Springer, New York, 2012.
	
	\bibitem{Cooper}
	J.N. Cooper, Constraints on Brouwer’s Laplacian spectrum conjecture, \emph{Linear Algebra Appl.},  \textbf{615} (2021) 11--27.
	
	\bibitem{DCve}D. Cvetkovi\'{c}, P. Rowlinson, S. Simi\'{c}, An Introduction to the Theory of Graph Spectra, Cambridge University Press, Cambridge, 2010.
	
	\bibitem{DuZhou}Z. Du, B. Zhou, Upper bounds for the sum of Laplacian eigenvalues of graphs, \emph{Linear Algebra Appl.},  \textbf{436} (2012) 3672--3683.
	
	\bibitem{KF}K. Fan, On a theorem of wely concerning eigenvalues of linear transformations, \emph{Proc. Natl. Acad. Sci.}, \textbf{35}(1949), 652--655.
	
	\bibitem{MF}M. Fiedler, Special Matrices and their Applications in Numerical Mathematics, Martinus Nijhoff, Dordrecht, 1986.
	
	\bibitem{Fritscher}E. Fritscher, C. Hoppen, I. Rocha, V. Trevisan, On the sum of the Laplacian eigenvalues of a tree, \emph{Linear Algebra and its Applications.}, \textbf{435}(2011), 371--399.
	
	\bibitem{Haemers}W.H. Haemers, A. Mohammadian, B. Tayfeh-Rezaie, On the sum of Laplacian eigenvalues of graphs, \emph{Linear Algebra Appl.}, \textbf{432}(2010), 2214--2221.
	
	\bibitem{LG}W. Li, J. Guo, On the full Brouwer's Laplacian spectrum conjecture, \emph{Discrete Mathematics.}, \textbf{345}(2022)113078.
	
	\bibitem{Mayank} Mayank, On Variants of the Grone-Merris Conjecture, Master’s thesis, Eindhoven University of Technology, 2010.
	
	\bibitem{Oliveira}C.S. Oliveira, L.S. Lima, P. Rama, P. Carvalho, Extremal graphs for the sum of the two largest signless Laplacian eigenvalues, \emph{Electronic Journal of Linear Algebra.}, \textbf{30}(2015), 605--612.
	
	\bibitem{WangHuang}S. Wang, Y. Huang, B. Liu, On a conjecture for the sum of Laplacian eigenvalues, \emph{Math. Comput. Model.}, \textbf{56} (2012) 60--68.
	
	\bibitem{WH}B. Wu, Y. Lou, C. He, Signless Laplacian and normalized Laplacian on the h-join operation of graphs, \emph{Discrete Mathematics Algorithms Appl.}, \textbf{6}(2014)1450046.
	
	\bibitem{Zheng}Y. Zheng, A note on the sum of the two largest signless Laplacian eigenvalues, \emph{Ars Combin.}, \textbf{148}(2020), 183--191.
	
	
	
	
	
\end{thebibliography}
\end{document}



 






\title{{ Corrigendum to ''On the sum of the first two largest signless Laplacian eigenvalues of a graph'' 
 
 [Discrete Math. 347 (2024) 114035]}
 \footnote{ DOI of original article: https://doi.org/10.1016/j.disc.2024.114035\\
            Corresponding author:\ Chang-Xiang He.\\
            Corresponding author Email:\ cxhe@usst.edu.cn (C.-X. He).}}

\author[1]{Zhi-Bin Du}		
\author[2]{Zi-Ming Zhou}
\author[3]{Hai-Ying Shan}
\author[2]{Chang-Xiang He}

\affil[1]{School of Software, South China Normal University, Foshan, Guangdong 528225,  P. R. China}
\affil[2]{College of Science, University of Shanghai for Science and Technology, Shanghai, P. R. China.}
\affil[3]{School of   Mathematical Sciences, Tongji University, Shanghai, P. R. China.}

  \renewcommand\Authands{ and }

 \date{}

 \maketitle

The purpose of this corrigendum is to rectify an error in our recent publication \cite{ZHS}, in which we wrongly claimed that $f(G_{t_1,t_2,t_3})>f(G_{t_1,t_2+1,t_3-1})$ with $t_2 \ge t_3$ aiming to $G_1 \cong G_{t_1,t_2,t_3}$, as well as $G_{11}$ and $G_{14}$. This error is attributed to the fact that
$$
P_1-P_2 = P(\Delta_2(Q^\pi(G_{t_1,t_2,t_3})),x+e(G_{t_1,t_2,t_3})+3)-P(\Delta_2(Q^\pi(G_{t_1,t_2+1,t_3-1})),x+e(G_{t_1,t_2+1,t_3-1})+3)
$$
is not always positive. In this corrigendum, we will give supplements and corrections to our previous proof. The notations not defined here can refer to \cite{ZHS}.

We have proved in \cite{ZHS} that
$\frac{1.3}{n}<f(K^+_{1,n-1})<\frac{1.5}{n}$, and $f(G)>\frac{1.5}{n}$ when $G \not\cong G_1, G_{11}, G_{14}$. In the following, we will show that $f(G)>\frac{1.5}{n}$ for $G \cong G_1$ with $(t_1, t_2, t_3) \ne (1,n -3, 0)$, or $G \cong G_{11}, G_{14}$. Then it results in $f(G) \ge f(K^+_{1,n-1})$ with equality if and only if $G \cong K^+_{1,n-1}$, completing the proof of \cite[Theorem 1.1]{ZHS}. It is worth mentioning that $K^+_{1,n-1}$ is the same as $G_{t_1,t_2,t_3}$ with $(t_1, t_2, t_3) = (1,n -3, 0)$.

\begin{figure}[htbp!]
\centering
\includegraphics[width=12cm,page={1}]{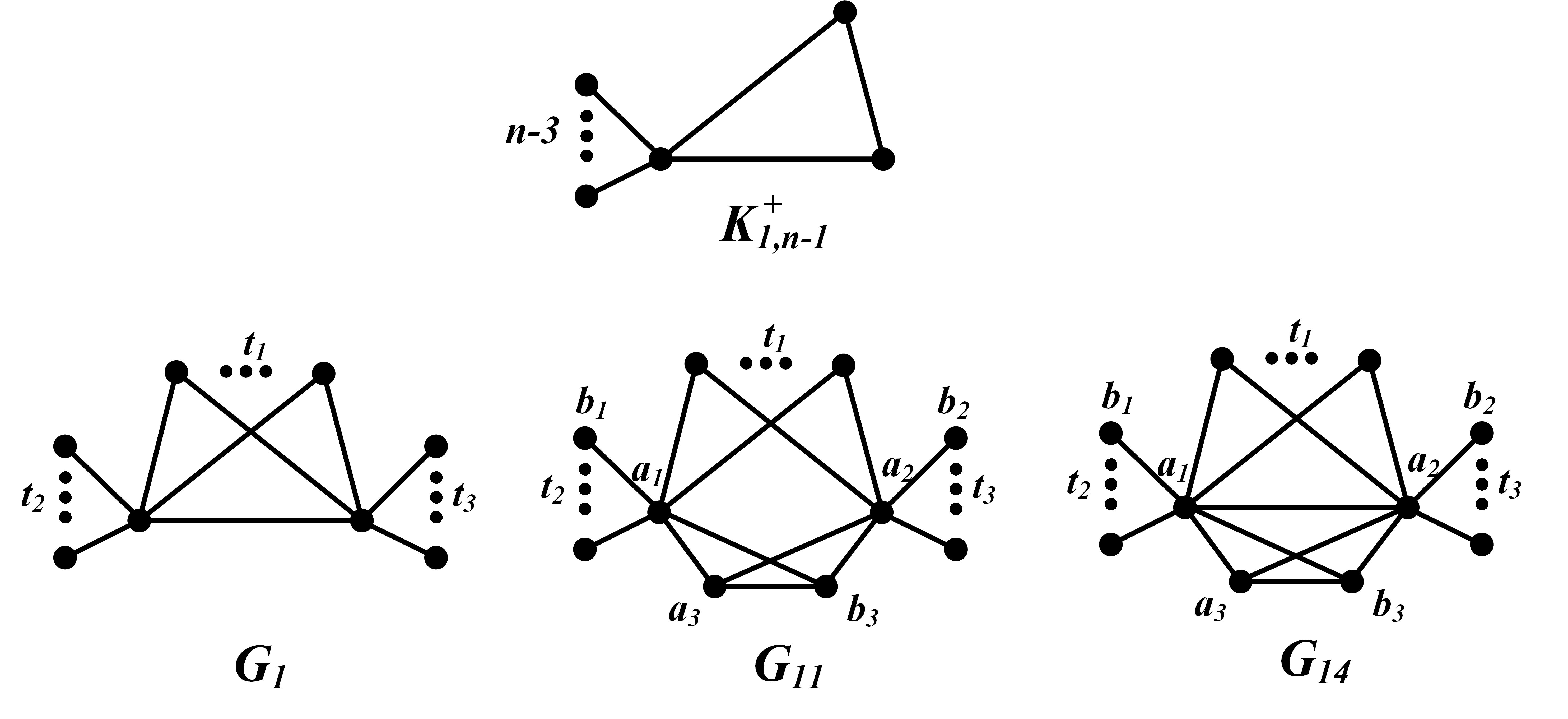}
\caption{The graphs $K^+_{1,n-1}$, $G_1$, $G_{11}$ and $G_{14}$}
\label{G1G11G14}
\end{figure}

\begin{case}
$G\cong G_1$ with $(t_1, t_2, t_3) \ne (1,n -3, 0)$.
\end{case}

It is easy to get that the characteristic polynomial of $Q(G_{t_1,t_2,t_3})$ is equal to
\begin{align*}
P_Q(G_{t_1,t_2,t_3},x)=P(Q^\pi(G_{t_1,t_2,t_3}),x)(x-2)^{t_1-1}(x-1)^{t_2+t_3-2}.
\end{align*}
On one hand, since $q_2>2$ by Lemma 2.9, $P(Q^\pi(G_{t_1,t_2,t_3}),x)$ should contain the largest two roots of $P_Q(G_{t_1,t_2,t_3},x)$. On the other hand, one can see that all the coefficients of
$$P\left(\Delta_2(Q^\pi(G_{t_1,t_2,t_3})),x+e(G)+3-\frac{1.5}{n}\right)
$$
on $x$ are positive (since the denominator is positive, we need only to consider the numerator after the reduction of fractions to a common denominator, the details are deferred in \cite{P}),
which infers that the largest eigenvalue of $\Delta_2(Q^\pi(G_{t_1,t_2,t_3}))$ is less than $e(G)+3-\frac{1.5}{n}$.
Now by Corollary 2.1, we have $f(G)>\frac{1.5}{n}$.

\begin{case}
$G\cong G_{11}$ or $G_{14}$.
\end{case}

We have checked $f(G)>\frac{1.5}{n}$ holds for all the graphs with at most $11$ vertices by computational checking.  Assuming  the graph under consideration has at least $12$ vertices  next.
Clearly, $d(a_3)+d(b_3)=6$.

According to the symmetry of the structure of $G_{11}$ ($G_{14}$, respectively),
without loss of generality,
if $t_1+t_2>2$, then $G_{11}$ ($G_{14}$, respectively) contains two disjoint stars $K_{1,5}$ with centers $a_1$ and $a_2$.
Further, by Lemma 2.9, $q_2(G-e_3)\geq 6$ ($e_3=a_3b_3$).
Using Lemmas 3.3 and 2.10, we find that 
$$
S_2(G)<S_2(G-e_3)+1<\left(e(G-e_3)+3-\frac{1.5}{n}\right)+1=e(G)+3-\frac{1.5}{n},
$$
equivalently $f(G)>\frac{1.5}{n}$.

Since $|N(M)|=2$, $t_1+t_2>1$, we are remaining  to consider $t_1+t_2=2$, that is to say, $(t_1, t_2) = (0,2), (1,1)$ or $(2,0)$. For each of such cases, using the $2^{nd}$ additive compound matrix approach, we have $f(G)>0.4>\frac{1.5}{n}$, the details can refer to \cite{ZHS}.

The authors would like to apologise for any inconvenience caused.